\documentclass[a4paper,11pt]{amsart}
\usepackage{graphicx}
\usepackage{mathptmx}
\usepackage{amsmath}
\usepackage{amssymb}
\usepackage{enumitem}
\usepackage{xcolor}
\usepackage{xparse}
\NewDocumentCommand{\eulerian}{omm}
 {%
  \genfrac<>{0pt}{}{#2}{#3}%
  \IfValueT{#1}{_{\!#1}}%
 }

 \newmuskip\pFqmuskip

\newcommand*\pFq[6][8]{%
  \begingroup 
  \pFqmuskip=#1mu\relax
  \mathchardef\normalcomma=\mathcode`,
  \mathcode`\,=\string"8000
  \begingroup\lccode`\~=`\,
  \lowercase{\endgroup\let~}\pFqcomma
  {}_{#2}F_{#3}{\left(\genfrac..{0pt}{}{#4}{#5}\bigg|#6\right)}%
  \endgroup
}
\newcommand{\pFqcomma}{{\normalcomma}\mskip\pFqmuskip}

\newtheorem{theorem}{Theorem}

\begin{document}

\title[On Some Summation Formulas]{On Some Summation Formulas}

\author{Taekyun  Kim}
\address{Department of Mathematics, Kwangwoon University, Seoul 139-701, Republic of Korea}
\email{tkkim@kw.ac.kr}

\author{DAE SAN KIM}
\address{Department of Mathematics, Sogang University, Seoul 121-742, Republic of Korea}
\email{dskim@sogang.ac.kr}

\author{Hyunseok Lee}
\address{Department of Mathematics, Kwangwoon University, Seoul 139-701, Republic of Korea}
\email{luciasconstant@kw.ac.kr}

\author{Jongkyum Kwon}
\address{Department of Mathematics Education, Gyeongsang National University, Jinju 52828, Republic of Korea}
\email{mathkjk26@gnu.ac.kr}


\maketitle

\begin{abstract}
In this note, we derive a finite summation formula and an infinite summation formula involving Harmonic numbers of order up to some order by means of several definite integrals.
\end{abstract}

\section{Introduction}
For $s\in\mathbb{C}$ with $\mathrm{Re}(s)>0$, the gamma function is defined by
\begin{equation}
\Gamma(s)=\int_{0}^{\infty}e^{-t}t^{s-1}dt,\quad (\mathrm{see}\ [5]).\label{1}
\end{equation}
From \eqref{1}, we note that $\Gamma(s+1)=s\Gamma(s)$ and $\Gamma(n+1)=n!$, $(n\in\mathbb{N})$, (see [1,2,3]). \par
As is well known, the Beta function is defined for $\mathrm{Re}\,\alpha>0,\, \mathrm{Re}\,\beta>0$ by
\begin{equation}
B(\alpha,\beta)=\int_{0}^{1}t^{\alpha-1}(1-t)^{\beta-1}dt\label{2},\quad (\mathrm{see}\ [1,5]).
\end{equation}
From \eqref{2} we note that
\begin{equation}
B(\alpha,\beta)=\frac{\Gamma(\alpha)\Gamma(\beta)}{\Gamma(\alpha+\beta)},\quad (\mathrm{see}\ [5]).\label{3}
\end{equation}
The Harmonic numbers are defined by
\begin{equation}
H_{n}=1+\frac{1}{2}+\frac{1}{3}+\cdots+\frac{1}{n}, \quad (n\ge 1),
\end{equation}
and more generally, for any $r\in\mathbb{N}$, the Harmonic numbers of order $r$ are given by
\begin{displaymath}
H_{n}^{(r)}=1+\frac{1}{2^{r}}+\frac{1}{3^{r}}+\cdots+\frac{1}{n^{r}},\quad (n\ge 1).
\end{displaymath} \par

In this note, we derive a finite summation formula in \eqref{A} and an infinite summation formula \eqref{B}, where $f_{n,r}(0)$ is determined by the recurrence relation in \eqref{C} and given by Harmonic numbers of order $j$, for $j=1,2,\dots,r$. Our results are illustrated for $r=3$ and $r=4$ in the Example below. It is amusing that the infinite sum in \eqref{B} involving Harmonic numbers of order $\le r$ boils down to $(-1)^{r}(r+1)!$. They are derived from several definite integrals in an elementary way. After considering the two summation formulas in \eqref{A} and \eqref{B} for $r=1$ and $r=2$, we will show the results in Theorem 1.

\begin{theorem}
Let $n, r$ be positive integers. Then we have
\begin{align}
&(-1)^{r}r!\sum_{k=0}^{n}\binom{n}{k}(-1)^{k}\frac{1}{(2k+1)^{r+1}}=\frac{2^{2n}}{(2n+1)\binom{2n}{n}}f_{n,r}(0), \label{A}\\
&\sum_{n=1}^{\infty}\frac{2^{2n-1}}{n(2n+1)\binom{2n}{n}}f_{n,r}(0)=(-1)^{r}(r+1)!, \label{B}
\end{align}
where $f_{n,s}(x)$ are determined by the recurrence relation
\begin{align}
f_{n,s+1}(x)=\beta_{n}(x)f_{n,s}(x)+\frac{d}{dx}f_{n,s}(x), \,\,(s \ge 1),\,\, f_{n,1}(x)=\beta_{n}(x)=-\sum_{k=0}^{n}(x+2k+1)^{-1}, \label{C}
\end{align}
so that $f_{n,s}(x)$ is a polynomial in $x$ involving $\beta_{n}(x), \beta_{n}^{(1)}(x),\cdots,\beta_{n}^{(s-1)}(x)$, with
\begin{align}
\beta_{n}^{(j)}(0)=(-1)^{j-1}j!\big(H_{2n+1}^{(j+1)}-\frac{1}{2^{j+1}}H_{n}^{(j+1)}\big),\quad(j \ge 0).\label{D}
\end{align}
\end{theorem}

\noindent{\bf{Example:}}\\
(a)
\begin{align*}
&-3!\sum_{k=0}^{n}\binom{n}{k}(-1)^{k}\frac{1}{(2k+1)^{4}}=\frac{2^{2n}}{(2n+1)\binom{2n}{n}}\\
& \times \left\{-(H_{2n+1}-\frac{1}{2}H_{n})^{3}-3(H_{2n+1}-\frac{1}{2}H_{n})(H_{2n+1}^{(2)}-\frac{1}{4}H_{n}^{(2)})-2(H_{2n+1}^{(3)}-\frac{1}{8}H_{n}^{(3)}) \right\}, \\
&\sum_{n=1}^{\infty}\frac{2^{2n-1}}{n(2n+1)\binom{2n}{n}}\left\{-(H_{2n+1}-\frac{1}{2}H_{n})^{3}-3(H_{2n+1}-\frac{1}{2}H_{n})(H_{2n+1}^{(2)}-\frac{1}{4}H_{n}^{(2)})-2(H_{2n+1}^{(3)}-\frac{1}{8}H_{n}^{(3)}) \right\} \\
&=-4!.
\end{align*}
(b)
\begin{align*}
& 4!\sum_{k=0}^{n}\binom{n}{k}(-1)^{k}\frac{1}{(2k+1)^{5}}=\frac{2^{2n}}{(2n+1)\binom{2n}{n}}
\left\{(H_{2n+1}-\frac{1}{2}H_{n})^{4}+6(H_{2n+1}-\frac{1}{2}H_{n})^{2}(H_{2n+1}^{(2)}-\frac{1}{4}H_{n}^{(2)}) \right. \\
& \quad\quad\quad\quad\quad\quad\quad\quad \left. 3(H_{2n+1}^{(2)}-\frac{1}{4}H_{n}^{(2)})^{2}+8(H_{2n+1}-\frac{1}{2}H_{n})(H_{2n+1}^{(3)}-\frac{1}{8}H_{n}^{(3)})+6(H_{2n+1}^{(4)}-\frac{1}{16}H_{n}^{(4)}) \right\}, \\
& 5!=\sum_{n=1}^{\infty}\frac{2^{2n-1}}{n(2n+1)\binom{2n}{n}}\left\{(H_{2n+1}-\frac{1}{2}H_{n})^{4}+6(H_{2n+1}-\frac{1}{2}H_{n})^{2}(H_{2n+1}^{(2)}-\frac{1}{4}H_{n}^{(2)}) \right. \\
& \quad\quad\quad\quad\quad\quad\quad\quad \left. 3(H_{2n+1}^{(2)}-\frac{1}{4}H_{n}^{(2)})^{2}+8(H_{2n+1}-\frac{1}{2}H_{n})(H_{2n+1}^{(3)}-\frac{1}{8}H_{n}^{(3)})+6(H_{2n+1}^{(4)}-\frac{1}{16}H_{n}^{(4)}) \right\}.
\end{align*}

\section{Derivation of Summation Formulas}
First, we observe that
\begin{equation}
\begin{aligned}
\int_{0}^{1}(1-x^{2})^{n}dx&=\sum_{k=0}^{n}\binom{n}{k}(-1)^{k}\int_{0}^{1}x^{2k}dx \\
&=\sum_{k=0}^{n}\binom{n}{k}(-1)^{k}\frac{1}{2k+1}.
\end{aligned}	\label{5}
\end{equation}
On the other hand,
\begin{align}
\int_{0}^{1}(1-x^{2})^{n}dx&=\frac{1}{2}\int_{0}^{1}(1-y)^{n}y^{-\frac{1}{2}}dy=\frac{1}{2}\int_{0}^{1}(1-y)^{n+1-1}\cdot y^{\frac{1}{2}-1}dy \label{6}\\
&=\frac{1}{2}B\bigg(n+1,\frac{1}{2}\bigg)=\frac{1}{2}\frac{\Gamma(n+1)\Gamma\big(\frac{1}{2}\big)}{\big(n+\frac{1}{2}\big)\big(n-\frac{1}{2}\big)\big(n-\frac{3}{2}\big)\cdots\frac{1}{2}\Gamma\big(\frac{1}{2}\big)}\nonumber \\
&=\frac{2^{n}n!}{(2n+1)(2n-1)(2n-3)\cdots 1}=\frac{2^{2n}n!n!}{(2n+1)(2n)!}=\frac{2^{2n}}{(2n+1)\binom{2n}{n}}.\nonumber
\end{align}
Thus, from \eqref{5} and \eqref{6} we note that
\begin{equation}
\binom{2n}{n}\sum_{k=0}^{n}\binom{n}{k}\frac{(-1)^{k}}{2k+1}=\frac{2^{2n}}{2n+1}.\label{7}
\end{equation}
From \eqref{7}, we have
\begin{equation}
\sum_{n=1}^{\infty}\frac{2^{2n}}{n(2n+1)\binom{2n}{n}}=\sum_{n=1}^{\infty}\frac{1}{n}\int_{0}^{1}(1-x^{2})^{n}dx=-2\int_{0}^{1}\log x\ dx=2.\label{8}
\end{equation}
Let
\begin{equation}
\begin{aligned}
F(x)&=\int_{0}^{1}(1-t^{2})^{n}t^{x}dt=\sum_{k=0}^{n}\binom{n}{k}(-1)^{k}\int_{0}^{1}t^{2k+x}dt=\sum_{k=0}^{n}\binom{n}{k}(-1)^{k}\frac{1}{2k+x+1}.\label{9}
\end{aligned}
\end{equation}
Note that
\begin{equation}
\begin{aligned}
F^{\prime}(x)=\frac{d}{dx} \int_{0}^{1}(1-t^{2})^{n}t^{x}dt=\int_{0}^{1}(1-t^{2})^{n}t^{x}\log t dt=-\sum_{k=0}^{n}\binom{n}{k}(-1)^{k}\frac{1}{(2k+1+x)^{2}}.\label{10}
\end{aligned}
\end{equation}
On the other hand
\begin{align}
\int_{0}^{1}(1-t^{2})^{n}t^{x}dt&=\frac{1}{2}\int_{0}^{1}(1-y)^{n}y^{\frac{x+1}{2}-1}dy=\frac{1}{2}B\bigg(n+1,\frac{x+1}{2}\bigg)	\label{11}\\
&=\frac{1}{2}\frac{\Gamma(n+1)\Gamma\big(\frac{x+1}{2}\big)}{\Gamma\big(n+1+\frac{x+1}{2}\big)}=\frac{1}{2}\frac{n!\Gamma\big(\frac{x+1}{2}\big)}{\big(n+\frac{x+1}{2}\big)\big(n+\frac{x-1}{2}\big)\big(n+\frac{x-3}{2}\big)\cdots\big(\frac{x+1}{2}\big)\Gamma\big(\frac{x+1}{2}\big)}\nonumber \\
&=\frac{n!2^{n}}{(2n+x+1)(2n+x-1)(2n+x-3)\cdots(x+1)}=\frac{n!2^{n}}{\prod_{k=0}^{n}(x+2k+1)}.\nonumber
\end{align}
Thus, we have
\begin{equation}
\begin{aligned}
F^{\prime}(x)&=\frac{d}{dx}F(x)=\frac{d}{dx}\int_{0}^{1}(1-t^{2})^{n}t^{x}dx=\frac{d}{dx}\bigg(\frac{n!2^{n}}{\prod_{k=0}^{n}(x+2k+1)}\bigg)\\
&=-\bigg(\frac{n!2^{n}}{\prod_{k=0}^{n}(x+2k+1)}\bigg)\bigg(\sum_{k=0}^{n}\frac{1}{x+2k+1}\bigg).
\end{aligned}	\label{12}
\end{equation}
Thus, by \eqref{12}, we get
\begin{equation}
F^{\prime}(0)=-\bigg(\frac{n!n!2^{2n}}{(2n+1)(2n)!}\bigg)\bigg(\sum_{k=0}^{n}\frac{1}{2k+1}\bigg)=-\frac{2^{2n}}{(2n+1)\binom{2n}{n}}\bigg(H_{2n+1}-\frac{1}{2}H_{n}\bigg).\label{13}	
\end{equation}
From \eqref{10} and \eqref{13}, we note that
\begin{align}
\binom{2n}{n}\sum_{k=0}^{n}\binom{n}{k}(-1)^{k}\frac{1}{(2k+1)^{2}}=\frac{2^{2n}}{2n+1}\bigg(H_{2n+1}-\frac{1}{2}H_{n}\bigg).\label{14}
\end{align}
By \eqref{14}, we have
\begin{align}
\sum_{n=1}^{\infty}\frac{2^{2n}}{n(2n+1)\binom{2n}{n}}\big(H_{2n+1}-\frac{1}{2}H_{n}\big)&=\sum_{n=1}^{\infty}\frac{1}{n}\sum_{k=0}^{n}\binom{n}{k}\frac{(-1)^{k}}{(2k+1)^{2}}\label{15}\\
&=\sum_{n=1}^{\infty}\frac{1}{n}\int_{0}^{1}\int_{0}^{1}(1-x^{2}y^{2})^{n}dxdy\nonumber\\
&=-2\int_{0}^{1}\int_{0}^{1}(\log x+\log y)dxdy=4.\nonumber
\end{align}
Thus, we have
\begin{align}
\sum_{n=1}^{\infty}\frac{2^{2n}}{n(2n+1)\binom{2n}{n}}\bigg(H_{2n+1}-\frac{1}{2}H_{n}\bigg)=4.\label{16}
\end{align}
From \eqref{9}, we note that
\begin{align}
F^{\prime\prime}(x)=\frac{d^{2}}{dx^{2}}F(x)&=\frac{d^{2}}{dx^{2}}\int_{0}^{1}(1-t^{2})^{n}t^{x}dt\label{17} \\
&=\int_{0}^{1}(1-t^{2})^{n}\big(\log t\big)^{2}t^{x}dt\nonumber \\
&=2!\sum_{k=0}^{n}\binom{n}{k}(-1)^{k}\frac{1}{(2k+1+x)^{3}}.\nonumber	
\end{align}
Hence, by \eqref{17}, we get
\begin{equation}
F^{\prime\prime}(x)=\int_{0}^{1}(1-t^{2})^{n}\big(\log t\big)^{2}t^{x}dt=2!\sum_{k=0}^{n}\binom{n}{k}(-1)^{k}\frac{1}{(2k+1+x)^{3}}.\label{18}
\end{equation}
From \eqref{11} and \eqref{12}, we have
\begin{equation}
\begin{aligned}
F^{\prime\prime}(x)&=\frac{d^{2}}{dx^{2}}\int_{0}^{1}(1-t^{2})^{n}t^{x}dt \\
&=\frac{n!2^{n}}{(2n+1+x)(2n+x-1)\cdots(x+1)}\bigg\{\bigg(\sum_{k=0}^{n}\frac{1}{2k+x+1}\bigg)^{2}+\sum_{k=0}^{n}\frac{1}{(2k+x+1)^{2}}\bigg\}.
\end{aligned}\label{19}
\end{equation}
Thus, we note that
\begin{equation}
\begin{aligned}
F^{\prime\prime}(0)& =\frac{n!2^{n}}{(2n+1)(2n-1)(2n-3)\cdots 1}\bigg\{\bigg(\sum_{k=0}^{n}\frac{1}{2k+1}\bigg)^{2}+\sum_{k=0}^{n}\frac{1}{(2k+1)^{2}}\bigg\}\\
&=\frac{n!n!2^{2n}}{(2n+1)(2n)!}\bigg\{\big(H_{2n+1}-\frac{1}{2}H_{n}\big)^{2}+\bigg(H_{2n+1}^{(2)}-\frac{1}{4}H_{n}^{(2)}\bigg)\bigg\}.
\end{aligned}\label{20}
\end{equation}
Therefore, by \eqref{17} and \eqref{20}, we get
\begin{align}
2!\binom{2n}{n}\sum_{k=0}^{n}\binom{n}{k}\frac{(-1)^{k}}{(2k+1)^{3}}=\frac{2^{2n}}{(2n+1)}\bigg\{\bigg(H_{2n+1}-\frac{1}{2}H_{n}\bigg)^{2}+\bigg(H_{2n+1}^{(2)}-\frac{1}{4}H_{n}^{(2)}\bigg) \bigg\}.\label{21}
\end{align}
Note that
\begin{equation}
\begin{aligned}
\sum_{n=1}^{\infty}\frac{1}{n}\bigg\{\sum_{k=0}^{n}\binom{n}{k}\frac{(-1)^{k}}{(2k+1)^{3}}\bigg\}&=\sum_{n=1}^{\infty}\frac{1}{n}\int_{0}^{1}\int_{0}^{1}\int_{0}^{1}(1-x^{2}y^{2}z^{2})^{n}dxdydz\\
&=-2\int_{0}^{1}\int_{0}^{1}\int_{0}^{1}(\log x+\log y+\log z)dxdydz=6.
\end{aligned}\label{22}
\end{equation}
From \eqref{21} and \eqref{22}, we have
\begin{align}
\sum_{n=1}^{\infty}\frac{2^{2n-1}}{n(2n+1)\binom{2n}{n}}\bigg\{\bigg(H_{2n+1}-\frac{1}{2}H_{n}\bigg)^{2}+ \bigg(H_{2n+1}^{(2)}-\frac{1}{4}H_{n}^{(2)}\bigg)^{2} \bigg\}	=6. \label{23}
\end{align}
Now, we begin to prove Theorem 1. Let $F(x)=\int_{0}^{1}(1-t^2)^{n}t^{x} dt$ be as in \eqref{9}.\\
Let $r$ be a positive integer. Then repeated integrating by parts gives us
\begin{align}
F^{(r)}(x)&=\bigg(\frac{d}{dx}\bigg)^{r}\int_{0}^{1}(1-t^2)^{n}t^{x} dt =\int_{0}^{1}(1-t^2)^{n}(\log t)^{r}t^{x} dt \nonumber \\
&=\sum_{k=0}^{n}\binom{n}{k}(-1)^{k}\int_{0}^{1}(\log t)^{r}t^{x+2k} dt \nonumber\\
&=\sum_{k=0}^{n}\binom{n}{k}(-1)^{k}\frac{-r}{x+2k+1}\int_{0}^{1}(\log t)^{r-1}t^{x+2k} dt \label{24}\\
&=\sum_{k=0}^{n}\binom{n}{k}(-1)^{k}\frac{-r}{x+2k+1}\cdots\frac{-1}{x+2k+1}\int_{0}^{1}t^{x+2k} dt \nonumber\\
&=(-1)^{r}r!\sum_{k=0}^{n}\binom{n}{k}(-1)^{k}\frac{1}{(x+2k+1)^{r+1}}.\nonumber
\end{align}
As we saw in \eqref{11} and \eqref{12}, $F(x)$ is alternatively expressed by
\begin{align}
F(x)=\int_{0}^{1}(1-t^2)^{n}t^{x} dt=\frac{n!2^{n}}{\prod_{k=0}^{n}(x+2k+1)},\,\,\,F(0)=\frac{2^{2n}}{(2n+1)\binom{2n}{n}},\label{25}
\end{align}
and its derivative is given by
\begin{align}
F^{(1)}(x)=F(x)\beta_{n}(x),\quad \beta_{n}=\beta_{n}(x)=-\sum_{k=0}^{n}(x+2k+1)^{-1}.\label{26}
\end{align}
Repeated differentiations give us
\begin{align}
&F^{(2)}(x)=F(x)(\beta_{n}^{2}+\beta_{n}^{(1)}),\,\,F^{(3)}(x)=F(x)(\beta_{n}^{3}+3\beta_{n}\beta_{n}^{(1)}+\beta_{n}^{(2)}), \nonumber\\
&F^{(4)}(x)=F(x)(\beta_{n}^{4}+6\beta_{n}^{2}\beta_{n}^{(1)}+3(\beta_{n}^{(1)})^{2}+4\beta_{n}\beta_{n}^{(2)}+\beta_{n}^{(3)}),
\label{27}\\
&F^{(5)}(x)=F(x)\left(\beta_{n}^{5}+10\beta_{n}^{3}\beta_{n}^{(1)}+10\beta_{n}^{2}\beta^{(2)}+15\beta(\beta_{n}^{(1)})^{2} \right.\nonumber \\
& \left.\quad\quad\quad\quad\quad\quad\quad+5\beta_{n}\beta_{n}^{(3)}
+10\beta_{n}^{(1)}\beta_{n}^{(2)}+\beta_{n}^{(4)}\right),\dots.\nonumber
\end{align}
In general, we let $F^{(s)}(x)=F(x)f_{n,s}(x)$, for $s \ge 1$. Then further differentiation gives us
\begin{align*}
F^{(s+1)}(x)=F(x)(\beta_{n}(x)f_{n,s}(x)+\frac{d}{dx}f_{n,s}(x)).
\end{align*}
Thus we get the recurrence relation for $\left\{f_{n,s}(x)\right\}_{s=1}^{\infty}$:
\begin{align}
f_{n,s+1}(x)=\beta_{n}(x)f_{n,s}(x)+\frac{d}{dx}f_{n,s}(x), \,\,(s \ge 1),\,\, f_{n,1}(x)=\beta_{n}(x).\label{28}
\end{align}
Now, from \eqref{24} and \eqref{25}, we obtain
\begin{align}
(-1)^{r}r!\sum_{k=0}^{n}\binom{n}{k}(-1)^{k}\frac{1}{(x+2k+1)^{r+1}}=\frac{n!2^{n}}{\prod_{k=0}^{n}(x+2k+1)}f_{n,r}(x). \label{29}
\end{align}
Letting $x=0$ in \eqref{29}, we get
\begin{align}
\sum_{k=0}^{n}\binom{n}{k}(-1)^{k}\frac{1}{(2k+1)^{r+1}}=\frac{(-1)^{r}}{r!}\frac{2^{2n}}{(2n+1)\binom{2n}{n}}f_{n,r}(0).\label{30}
\end{align}
Multiplying both sides of \eqref{30} by $\frac{1}{n}$ and summing over $n$, we get
\begin{align}
\frac{(-1)^{r}}{r!}\sum_{n=1}^{\infty}\frac{2^{2n}}{n(2n+1)\binom{2n}{n}}f_{n,r}(0)&=\sum_{n=1}^{\infty}\frac{1}{n}\sum_{k=0}^{n}\binom{n}{k}(-1)^{k}\frac{1}{(2k+1)^{r+1}} \nonumber \\
&=\sum_{n=1}^{\infty}\frac{1}{n}\int_{0}^{1} \cdots \int_{0}^{1}(1-x_{1}^{2}x_{2}^{2}\cdots x_{r+1}^{2})^{n} dx_{1}dx_{2}\cdots dx_{r+1}\label{31}\\
&=-2\int_{0}^{1} \cdots \int_{0}^{1}(\log x_1 +\cdots+\log x_{r+1})dx_{1} \cdots dx_{r+1}\nonumber\\
&=2(r+1).\nonumber
\end{align}
Finally, we note that $f_{n,s}(x)$ is a polynomial in $x$ involving $\beta_{n}(x), \beta_{n}^{(1)}(x),\cdots,\beta_{n}^{(s-1)}(x)$, as we can see from the recurrence relation in \eqref{28}. From \eqref{26}, we see that
\begin{align}
\beta_{n}^{(j)}(0)=(-1)^{j-1}j!\sum_{k=0}^{n}(2k+1)^{-(j+1)}=(-1)^{j-1}j!\big(H_{2n+1}^{(j+1)}-\frac{1}{2^{j+1}}H_{n}^{(j+1)}\big).\label{32}
\end{align}
Combining \eqref{28} and \eqref{30}-\eqref{32} altogether, we obtain our main result in Theorem 1.

\end{document}